# Using Edge-induced and Vertex-induced Subhypergraph Polynomials

Yohannes Tadesse

**Abstract**

For a hypergraph $\mathcal{H}$, we consider the edge-induced and vertex-induced subhypergraph polynomials and study their relation. We use this relation to prove that both polynomials are reconstructible, and to prove a theorem relating the Hilbert series of the Stanley-Reisner ring of the independent complex of $\mathcal{H}$ and the edge-induced subhypergraph polynomial. We also consider reconstruction of some algebraic invariants of $\mathcal{H}$.

## 1 Introduction

To every hypergraph $\mathcal{H}$ one can associate several subhypergraph enumerating polynomials. In this note we consider two of these polynomials: the vertex-induced subhypergraph polynomial $P_\mathcal{H}(x, y)$ enumerating vertex-induced subhypergraphs of $\mathcal{H}$, and the edge-induced subhypergraph polynomial $S_\mathcal{H}(x, y)$. Precise definitions will be given in §2. These and several other polynomials were extensively studied for graphs, see [1, 4, 5, 8] and their citations. The notion has been naturally generalized to hypergraphs, see White [14].

L. Borzacchini, et al. [5] studied the relation between these and other subgraph enumerating polynomials. He earlier proved that both are reconstructible, i.e. they can be derived from the subgraph enumerating polynomials of vertex-deleted subgraphs, see [3, 4]. A. Goodarzi [9] used $S_\mathcal{H}(x, y)$ to compute the Hilbert series of the Stanley-Resiner ring of the independent complex of $\mathcal{H}$. More precisely, if $R$ is such a ring, then its Hilbert series $H_R(t)$ is given by

$$(1.1) \qquad H_R(t) = \frac{S_\mathcal{H}(t, -1)}{(1-t)^n}$$

where $n$ is the number of vertices in $\mathcal{H}$.

In section 2, we define the polynomials, and then prove that

$$S_\mathcal{H}(x, y) = (1-x)^n P_\mathcal{H}(\frac{x}{1-x}, 1+y).$$

In section 3, we use this relation to give a short and elementary proof of (1.1). One may compare our proof with the technical proof in [9]. In section 4, generalizing Borzacchini's results [3, 4], we prove that both polynomials are reconstructible for hypergraphs. We also prove the reconstruction problems of some algebraic invariants of the independent complex of $\mathcal{H}$, where their graph counter part is proven by Dalili, Faridi and Traves in [6]. That is, we consider reconstructibility of the Hilbert series, the $f$-vector, the (multi-)graded Betti numbers and some graded Betti tables of the independent complex of $\mathcal{H}$.

## 2 Preliminaries

A hypergraph is a pair $\mathcal{H} = (V, E)$ where $V$ is a set of elements called vertices and $E \subset 2^V$ is a set of distinct subsets of $V$ called edges such that for any two edges $\varepsilon_1, \varepsilon_2 \in E$, we have $\varepsilon_1 \subset \varepsilon_2 \Rightarrow \varepsilon_1 = \varepsilon_2$. A hypergraph $\mathcal{H}$ is called finite if the vertex set $V$ is finite. We say $\mathcal{H}$ is a $d$-hypergraph if $|\varepsilon| = d$ for each $\varepsilon \in E$, where $|\varepsilon|$ is the cardinality of $\varepsilon$. A graph is a 2-hypergraph. In this note we always consider finite hypergraphs.

Let $\mathcal{H} = (V, E)$ be hypergraph, $W \subset V$ and $L \subset E$. We say that $\mathcal{L} = (W, L)$ an *edge-induced subhypergraph* of $\mathcal{H}$ if $W = \cup_{\varepsilon \in L} \varepsilon$. We say that $\mathcal{H}_W = (W, L)$ is *vertex-induced subhypergraph* if $L$ is the largest subset of $E$ such that $L \subset 2^W$.

Let $\mathcal{H}$ be a hypergraph. The *edge-induced subhypergraph polynomial* $S_\mathcal{H}(x, y)$ is defined by

$$(2.1) \qquad S_\mathcal{H}(x, y) = \sum_{i,j} \theta_{ij} x^i y^j,$$

where $\theta_{00} = 1$ and for $i, j \geq 0$, $\theta_{ij}$ is the number of edge induced subhypergraphs of $\mathcal{H}$ with $i$ vertices and $j$ edges. Similarly, the *vertex-induced subhypergraph polynomial* $P_\mathcal{H}(x, y)$ of $\mathcal{H}$ is defined by

$$(2.2) \qquad P_\mathcal{H}(x, y) = \sum_{i,j} \beta_{ij} x^i y^j,$$

where $\beta_{00} = 1$ and for $i, j \geq 0$, $\beta_{ij}$ is the number of vertex induced subhypergraphs of $\mathcal{H}$ with $i$ vertices and $j$ edges.

We recall some simple properties of these polynomials. In what follows, $F_\mathcal{H}(x, y)$ refers to any one of the two polynomials.

1. If the hypergraph has connected components $\mathcal{H}_1, \ldots, \mathcal{H}_m$, we have $F_\mathcal{H}(x, y) = F_{\mathcal{H}_1}(x, y) \cdots F_{\mathcal{H}_m}(x, y)$. We also have $F(0, y) = 1$. If $E = \emptyset$, then $F_\mathcal{H}(x, y) = (1 + x)^n$.

2. $\sum_j \beta_{ij} = \binom{n}{i}$ and $\sum_i \theta_{ij} = \binom{m}{j}$ where $m$ is the number of edges in $\mathcal{H}$.

3. $S_\mathcal{H}(x, 0)$ is a subgraph polynomial of the 0-subhypergraphs, i.e. isolated vertices. $P_\mathcal{H}(x, 0)$ the polynomial of the independent subsets, i.e. sets of vertices having no edges in common.

4. If $\mathcal{H} = K_n$ is the complete graph, then $P_\mathcal{H}(x, y) = \sum_{i=0}^n \binom{n}{j} x^i y^{\binom{j}{2}}$ and if $\mathcal{H}$ is a star with $m$ edges, then $S_\mathcal{H}(x, y) = \sum_{j=0}^m \binom{m}{j} x^{j+1} y^j$.

The following Proposition is a generalization of Borzacchini [3]. Even though he considered graphs, the proofs can easily be generalized to hypergraphs.

**Proposition 2.1.** *Let $\mathcal{H}$ be a hypergraph on $n$ vertices. Then*
$$S_\mathcal{H}(x,y) = (1-x)^n P_\mathcal{H}(\frac{x}{1-x}, 1+y)$$

*Proof.* To every vertex induced subhypergraph with $i$ vertices and $l$ edges there are $\binom{l}{j}$ hypergraphs with $i$ vertices and $j$ edges. Moreover, those obtained from different vertex induced subhypergraphs are different since they contain different vertex sets. On the other hand, to every edge induced subhypergraph with $l$ vertices and $j$ edges we can construct $\binom{n-l}{i-j}$ hypergraphs with $i$ vertices and $j$ edges. So

$$(2.3) \quad \sum_{l=0}^{i} \beta_{i,j+l} \binom{j+l}{j} = \sum_{l=0}^{i} \theta_{i-l,j} \binom{n-(i-l)}{l}.$$

Setting $r = j + l$ and $s = i - l$, substituting this in (2.3) and multiplying it by $x^i y^j$, we obtain:

$$\sum_{i,j} x^i y^j \left[\sum_{l=0}^{i} \beta_{i,j+l} \binom{j+l}{j}\right] = \sum_{i,j} x^i y^j \left[\sum_{l=0}^{i} \theta_{i-l,j} \binom{n-(i-l)}{l}\right].$$

$$\sum_{i,j} x^i y^j \left[\sum_{r} \beta_{ir} \binom{r}{j}\right] = \sum_{s,l,j} x^{s+l} y^j \left[\sum_{l=0}^{i} \theta_{sj} \binom{n-s}{l}\right].$$

$$\sum_{i,r} \beta_{ir} x^i \left[\sum_{j} \binom{r}{j} y^j\right] = \sum_{s,j} \theta_{sj} x^s y^j \left[\sum_{l} x^l \binom{n-s}{l}\right].$$

$$\sum_{i,r} \beta_{ir} x^i (1+y)^r = \sum_{s,j} \theta_{sj} x^s y^j (1+x)^{n-s}.$$

$$P_\mathcal{H}(x, y+1) = (1+x)^n \sum_{s,j} \theta_{sj} (\frac{x}{1+x})^s y^j.$$

$$P_\mathcal{H}(x, y+1) = (1+x)^n S_\mathcal{H}(\frac{x}{1+x}, y).$$

By change of variable, we obtain $S_\mathcal{H}(x,y) = (1-x)^n P_\mathcal{H}(\frac{x}{1-x}, 1+y)$. □

**Corollary 2.2.** *Let $\mathcal{H}$ be a hypergraph on $n$ vertices. Then*
$$P_\mathcal{H}(x,y) = (1+x)^n S_\mathcal{H}(\frac{x}{1+x}, y-1).$$

## 3 $P_\mathcal{H}(x,y)$ and $S_\mathcal{H}(x,y)$ in Algebra

A *simplicial complex* $\Delta$ on a vertex set $V = \{v_1, \ldots, v_n\}$ is a set of subsets of $V$, called faces or simplices such that $\{v_i\} \in \Delta$ for each $i$ and every subset of a face is itself a face. If $B \subset V$, the restriction of $\Delta$ to $B$ is a simplicial complex defined by $\Delta(B) = \{\delta \in \Delta \mid \delta \subset B\}$. The dimension of a face $\delta \in \Delta$ is $|\delta| - 1$. Let $f_i = f_i(\Delta)$ denote the number of faces of $\Delta$ of dimension $i$. Setting $f_{-1} = 1$, the sequence $f(\Delta) = (f_{-1}, f_0, f_1, \ldots, f_{d-1})$ is called the $f$-vector of $\Delta$.

Let $A = \mathbb{K}[x_1, \ldots, x_n]$ be a polynomial ring over a field $\mathbb{K}$ and $\Delta$ be a simplicial complex over $n$ vertices $V = \{v_1, \ldots, v_n\}$. The Stanley Reisner ideal of $\Delta$ is the ideal $I(\Delta) \subset A$ generated by those square free monomials $x_{i_1} \cdots x_{i_m}$ where $\{v_{i_1}, \ldots, v_{i_m}\} \notin \Delta$.

Let $\mathcal{H} = (V, E)$ be a hypergraph with $n$ vertices $V = \{v_1, \ldots, v_n\}$. An independent set of $\mathcal{H}$ is a subset $W \subset V$ such that $\varepsilon \not\subset W$ for all $\varepsilon \in E$. The collection of $\Delta_\mathcal{H}$ of independent sets forms a simplicial complex, called the *independent complex*. Thus the Stanley Resiner ideal of $\Delta_\mathcal{H}$ is the edge ideal of $\mathcal{H}$. More precisely, $I(\Delta_\mathcal{H}) = I(\mathcal{H}) \subset A$ is the ideal generated by the squarefree monomials $\prod_{x \in \varepsilon} x$ where $\varepsilon \in E$. Conversely, every squarefree monomial ideal $I \subset A$ can be associated with a hypergraph $\mathcal{H}_I = (V, E)$ where $V = \{v_1, \ldots, v_n\}$ and $\varepsilon \in E$ if $\prod_{x \in \varepsilon} x$ is in the minimal generating set of $I$. So one has $I(\Delta_{\mathcal{H}_I}) = I$. We have the following lemma.

**Lemma 3.1.** *Let $(f_0, f_1, \ldots, f_{d-1})$ be the $f$-vector of the independent complex of a hypergraph $\mathcal{H}$. Then $P_\mathcal{H}(t, 0) = \sum_{i=0}^{d} f_{i-1} t^i$.*

Let $R = \oplus_{i \in \mathbb{N}} R_n$ be a finitely generated graded $\mathbb{K}$-algebra, where $R_0 = \mathbb{K}$ is a field. The Hilbert series of $R$ is the generating function defined by $H_R(t) = \sum_{i \in \mathbb{N}} \dim_\mathbb{K}(R_i) t^i$. If $I \subset A$ is a monomial ideal, the Hilbert series of the monomial ring $R = A/I$ is the rational function $H_R(t) = \frac{\mathcal{K}(R,t)}{(1-t)^n}$ where $\mathcal{K}(R, t) \in \mathbb{Z}[t]$. P. Renteln [13], and also D. Ferrarello and R. Fröberg [7] used the subgraph induced polynomial $S_G(x, y)$ of a graph $G$ to compute the Hilbert series of the Stanley-Reisner ring $R$ of the independent complex of $G$, namely:

$$H_R(t) = \frac{S_G(t, -1)}{(1-t)^n}.$$

Recently A. Goodarzi [9] generalized it for any squarefree monomial ideal by using the combinatorial Alexander duality and Hochster's formula. Below is a very short and direct proof of this result.

**Theorem 3.2.** *Let $\mathcal{H}$ be a hypergraph on $n$ vertices, $I_\mathcal{H} \subset A = \mathbb{K}[x_1, \ldots, x_n]$ be its associated squarefree monomial ideal, and $R = A/I_\mathcal{H}$. Then*

$$H_R(t) = \frac{S_\mathcal{H}(t, -1)}{(1-t)^n}.$$

*Proof.* We know by Lemma 3.1 that $P_\mathcal{H}(t, 0) = \sum_{i=0}^{d} f_{i-1} x^i$ is the polynomial of the $f$-vectors of the independent complex of $\mathcal{H}$. It follows that by [12, Proposition 51.3] that $H_R(t) = P_\mathcal{H}(\frac{t}{1-t}, 0)$ and by Theorem 2.1 we have

$$S_\mathcal{H}(t, -1) = (1-t)^n P_\mathcal{H}(\frac{t}{1-t}, 0) = H_R(t)(1-t)^n.$$

□

**Remark 3.3.** Let $\mathcal{H}$ be a hypergraph and $R = A/I_\mathcal{H}$. It then follows by Lemma 3.1 and [12, Proposition 51.2] that $P_\mathcal{H}(t, 0)$ is the Hilbert polynomial of the exterior algebra $R/(x_1^2, \ldots, x_n^2)$.

# 4  $P_{\mathcal{H}}(x,y)$ and $S_{\mathcal{H}}(x,y)$ in reconstruction conjecture

For a graph $G = (V, E)$ on a vertex set $V = \{v_1, \ldots, v_n\}$, the deck of $G$ is the collection $\mathcal{D}(G) = \{G_1, \ldots, G_n\}$ where $G_l = G - v_l$, $v_l \in V$ is the vertex deleted subgraph of $G$. An element of $\mathcal{D}(G)$ is called a card. The long standing graph reconstruction conjecture posed by Kelly and Ulam says that every simple graph on $n \geq 3$ vertices is uniquely determined, up to isomorphism, by its deck. Numerous unsuccessful attempts have been made to prove the conjecture, nevertheless, a significant amount of work has been made. The reader may see Bondy [2] for a survey on the subject. Reconstruction of hypergraphs is defined similarly to graphs. Kocay [10] and Kocay and Lui [11] have constructed a family of non-reconstructible 3-hypergraphs.

**Remark 4.1.** Another obvious example of non-constructible hypergraphs are the 0-hypergraph containing no edges, and the $n$-hypergraph containing one edge with $n$-elements. So all the hypergraphs under consideration in this section are neither of these two.

In recent years questions has been asked if a graph invariant is reconstructible, that is, if it can be obtained from the its deck. Borzacchini in [3, 4] proved that both $S_G(x, y)$ and $P_G(x, y)$ are reconstructible. In fact, he proved that if $F_G(x, y)$ is any one of the subgraph polynomials and $F_{G_l}(x, y)$ is a subgraph polynomial of the card $G_l$, then

$$(4.1) \qquad nF_G(x,y) = x\frac{\partial F_G(x,y)}{\partial x} + \sum_{l=1}^{n} F_{G_l}(x,y).$$

It is natural to extend this reconstructibility question to hypergraphs. Below we obtain a similar result.

**Proposition 4.2.** *Let $\mathcal{H}$ be a hypergraph on $n \geq 3$ vertices. Then both $S_{\mathcal{H}}(x,y)$ and $P_{\mathcal{H}}(x,y)$ are reconstructible.*

*Proof.* We prove the proposition for $S_{\mathcal{H}}(x, y)$ since the other will follow by Proposition 2.1. Let $S_{\mathcal{H}}(x,y) = \sum_{ij} \theta_{ij} x^i y^j$ and $S_{\mathcal{H}_l}(x,y) = \sum_{ij} \theta_{ij}^{(l)} x^i y^j$ for $l = 1, \ldots, n$. By direct calculation we have

$$nS_{\mathcal{H}}(x,y) - x\frac{\partial(S_{\mathcal{H}}(x,y))}{\partial x} = n + \sum_{l=1}^{n}\sum_{ij}(n-j)\theta_{ij} x^i y^j.$$

Now if $j < n$, then any edge induced subhypergraph with $i$ vertices and $j$ edges is an edge induced subhypergraph for $n - j$ cards. It follows that $\sum_{l=1}^{n} \theta_{ij}^{(l)} = (n-j)\theta_{ij}$. Putting this in the equation and recalling that $n = \sum_{l=1}^{n} \theta_{00}^{(l)}$ we obtain

$$(4.2) \qquad nS_{\mathcal{H}}(x,y) = x\frac{\partial S_{\mathcal{H}}(x,y)}{\partial x} + \sum_{i=1}^{n} S_{\mathcal{H}_i}(x,y).$$

□

## 4.1 Hilbert series and $f$-vector

The authors in [6] studied reconstructibility of some algebraic invariants of the edge ideal of a graph $G$ such as the Krull dimension, the Hilbert series, and the graded Betti numbers $b_{i,j}$, where $j < n$. We extend these results to hypergraphs.

**Proposition 4.3.** *Let $\mathcal{H}$ be a hypergraph on $n \geq 3$ vertices. The Hilbert function of $R = A/I_\mathcal{H}$ is reconstructible.*

*Proof.* By Proposition 3.2 and (4.2) we have

$$nH_R(t) = \frac{nS_\mathcal{H}(t,-1)}{(1-t)^n} = \frac{t\frac{dS_\mathcal{H}(t,-1)}{dt}}{(1-t)^n} + \sum_{i-1}^{n} \frac{S_{\mathcal{H}_i}(t,-1)}{(1-t)^n}$$

$$= \frac{t}{(1-t)^n}\frac{dS_\mathcal{H}(t,-1)}{dt} + \sum_{i=1}^{n} \frac{H_{R_i}(t)}{1-t}.$$

Since $\frac{dH_R(t)}{dt} = \frac{d}{dt}\left(\frac{S_\mathcal{H}(t,-1)}{(1-t)^n}\right) = \frac{1}{t}\frac{t}{(1-t)^n}\frac{dS_\mathcal{H}(t,-1)}{dt} + \frac{n}{1-x}H_R(t)$, substituting this into the above, we obtain a first order ordinary linear differential equation

$$\frac{n}{1-t}H_R(t) = t\frac{dH_R(t)}{dt} - \frac{1}{1-t}\sum_{i=1}^{n}H_{R_i}(t).$$

□

**Proposition 4.4.** *Let $\mathcal{H}$ be a hypergraph on $n \geq 3$ vertices. The $f$-vector of $\Delta_\mathcal{H}$ is reconstructible.*

*Proof.* This, in fact, follows from Proposition 4.3, but we give an independent proof. Let $f(\Delta_\mathcal{H}) = (f_0, \ldots, f_{d-1})$. If $d < n$, by (4.2) when $F = P_\mathcal{H}$ and Lemma 3.1, we have $nf_{l-1} = if_{l-1} + \sum_{i=1}^{n} f_{i-1}^l$ for all $l \leq d$. If $d = n$, then $\mathcal{H}$ has no edges so $f_{d-1} = 1$. □

Let $\Delta_\mathcal{H}$ be the independent complex of a hypergraph $\mathcal{H}$. We can compute other invariants of $\Delta_\mathcal{H}$ from its $f$-vector $f(\Delta_\mathcal{H}) = (f_{-1}, f_0, \ldots, f_{d-1})$. Recall, for example, that the $h$-vector $h(\Delta_\mathcal{H}) = (h_0, \ldots, h_d)$ is defined by the formula $\sum_{i=0}^{d} f_{i-1}t^i(1-t)^{d-i} = \sum_{i=0}^{d} h_i t^i$. We can also obtain the multiplicity of the $R = A/I_\mathcal{H}$, namely $e(R) = f_{d-1}$. The following are consequences of Propositions 4.3 and 4.4.

**Corollary 4.5.** *Let $\mathcal{H}$ be a hypergraph on $n \geq 3$ vertices. The $h$-vector of $\Delta_\mathcal{H}$ is reconstructible.*

**Corollary 4.6.** *Let $\mathcal{H}$ be a hypergraph on $n \geq 3$ vertices. Then the Krull dimension and the multiplicity of $R = A/I_\mathcal{H}$ are reconstructible.*

## 4.2 Multi-graded Betti numbers

In this subsection we assume that $\operatorname{char} \mathbb{K} = 0$. Let $I \subset A = \mathbb{K}[x_1, \ldots, x_n]$ be a monomial ideal and consider the $\mathbb{Z}^n$-graded minimal free resolution of the $A$-module $R = A/I$:

$$\cdots \to \oplus_j A(-\mathbf{b})^{b_{i,\mathbf{b}}} \to \cdots \to \oplus_j A(-\mathbf{b})^{b_{2,\mathbf{b}}} \to \oplus_j A(-\mathbf{b})^{b_{1,\mathbf{b}}} \to A \to A/I \to 0$$

where $\mathbf{b} \in \mathbb{Z}^n$ and the modules $A(-\mathbf{b})$ are the shifts of $A$ to make the multi-graded differentials degree zero maps. The numbers $b_{i,\mathbf{b}}$ are *multi-graded Betti numbers* and $b_{ij} = \sum_{|\mathbf{b}|=j} b_{i,\mathbf{b}}$, where $|\mathbf{b}| = b_1 + \cdots + b_n$, are the *graded Betti numbers* of $R$. In particular, the $b_{in}$'s are the *extremal graded Betti numbers*. The importance of the assumption that char $\mathbb{K} = 0$ is that these numbers depend on the characteristic of the ground field, see eg. [12, Example 12.4]. If $I = I_{\mathcal{H}}$ is the edge ideal of a hypergraph $\mathcal{H}$, then each $\mathbf{b} \in \{0,1\}^n$, see for example [12, Corollary 26.10]. One can use graded Betti numbers to compute the Hilbert series of $R = A/I_{\mathcal{H}}$. So by Theorem 3.2, we have

$$(4.3) \qquad S_{\mathcal{H}}(t,-1) = \sum_{i=0}^{n} \sum_{j} (-1)^i b_{ij} t^j.$$

We generalize [6, Theorem 5.1] with a similar proof.

**Proposition 4.7.** *Let $\mathcal{H}$ be a hypergraph on with a vertex set $V = \{v_1, \ldots, v_n\}$ and $n \geq 3$. Then the multi-graded Betti numbers $b_{ij}$ of the Stanley Reisner ring $R = A/I_{\mathcal{H}}$ are reconstructible for all $j < n$.*

*Proof.* Let $\Delta = \Delta_{\mathcal{H}}$, $\Delta^{(l)} = \Delta_{\mathcal{H}_l}$, $\mathbf{b} \in \mathbb{Z}^n$, $b_{i,\mathbf{b}}$ be the multi-graded Betti numbers of $\Delta$, and $b_{i,\mathbf{b}}^{(l)}$ be the multi-graded Betti numbers of $\Delta^{(l)}$. By Hochester's formula, we have
$$b_{i,\mathbf{b}} = b_{i,B} = \tilde{h}_{j-i-1}(\Delta(B)),$$
where $B = \{v_i \in V \mid b_i \neq 0\}$ and $\tilde{h}_{j-i-1}(\Delta(B)) = \dim_{\mathbb{K}}(\tilde{H}_{j-i-1}(\Delta(B);\mathbb{K}))$ is the reduced simplicial homology of the subcomplex $\Delta(B)$. Since $\Delta(B) = \Delta^{(l)}(B)$ whenever $v_l \notin B$, it follows by Hochester's formula that $b_{i,\mathbf{b}} = \tilde{h}_{j-i-1}(\Delta^{(l)}(B)) = b_{i,\mathbf{b}}^{(l)}$. So the result holds. $\square$

**Corollary 4.8.** *Let $\mathcal{H}$ be a hypergraph with a vertex set $V = \{v_1, \ldots, v_n\}$ and $n \geq 3$. Then the graded Betti numbers $b_{ij}$ of the Stanley Reisner ring $R = A/I_{\mathcal{H}}$ are reconstructible for all $j < n$.*

*Proof.* $b_{ij} = \sum_{|\mathbf{b}|=j} b_{i,\mathbf{b}}$ and multi-graded Betti numbers are reconstructible. $\square$

Reconstruction of the extremal graded Betti numbers seems a bit hard to determine. We know that the coefficient of $t^n$ in $S_{\mathcal{H}}(t,-1)$ is the alternating sum $\sum_i (-1)^i b_{in}$. It follows that $b_{in}$ is reconstructible if there is only one $i$ such that $b_{in} \neq 0$. Fortunately, we have a good class of ideals with this property: for example, edge ideals of complements of chordal graphs, metroidal ideals, ideals with linear quotients and Cohen-Macaulay ideals. However, there are also edge ideals with more than one non-zero extremal graded Betti numbers [6, Example 5.3]. On the other hand, it is a useful invariant since it gives us information on many other invariants of $I_{\mathcal{H}}$. The following extends [6, Proposition 5.4] to hypergraphs.

**Proposition 4.9.** *Let $\mathcal{H}$ be a hypergraph on $n \geq 3$ vertices. If the graded top degree Betti numbers $b_{in}$ of $I_{\mathcal{H}}$ are reconstructible, then the depth, projective dimension and regularity of $I_{\mathcal{H}}$ are reconstructible.*

We investigate if the Betti table of $I_\mathcal{H}$ is reconstructible. Let $\mathcal{B} = (b_{ij})$ be the Betti table of $I_\mathcal{H}$ and $\mathcal{B}_l = (b_{ij}^{(l)})$ be the Betti table of $I_{\mathcal{H}_l}$. Then combining (4.2) and (4.3) and comparing the coefficients of $t^j$ we obtain

$$(n-j)\sum_i (-1)^i b_{ij} = \sum_i (-1)^i \sum_{l=1}^n b_{ij}^{(l)} \quad \text{for } j < n.$$

This equation shows it is difficult to determine each $b_{ij}$ only from the data $\{\mathcal{B}_l\}_{l=1}^n$ since anti-diagonals of $\mathcal{B}$ might contain more than one non-zero entry. We thus have the following which gives a partial answer to [6, Question 5.6].

**Proposition 4.10.** *Let $\mathcal{H}$ be a hypergraph on $n \geq 3$ vertices. If each anti-diagonal of the Betti table of $I_\mathcal{H}$ contains at most one non-zero entry, then the Betti table of $I_\mathcal{H}$ is reconstructible.*

In fact, in this case, we can compute the non-zero entries from the coefficients of $S_\mathcal{H}(x, y)$.

**Proposition 4.11.** *Let $\mathcal{H}$ be a hypergraph on $n \geq 3$ vertices and $S_\mathcal{H}(x, y) = \sum_{ij} \theta_{ij} x^i y^j$. Assume that each anti-diagonal of the Betti table contains at most one non-zero entry $b_0, b_1, \ldots, b_d$. Then $b_i = \sum_j \theta_{ij}$.*

*Proof.* Since $S_\mathcal{H}(t, -1) = \sum_{ij} \theta_{ij}(-1)^j t^i = \sum_{j=0}^n (-1)^j b_i t^i$, $b_i$ is the coefficient of $t^i$ in $S_\mathcal{H}(t, -1)$. □

*Yohannes Tadesse*
*Department of Mathematics*
*Uppsala University*
*Box: 480*
*SE 75206 Uppsala, Sweden*
*email: yohannes.tadesse@math.uu.se*